\documentstyle[10pt,amscd,xypic,amssymb,combelow,enumitem]{amsart}

\xyoption{all}
\CompileMatrices

\makeatletter
\@addtoreset{equation}{section}
\makeatother

\newtheorem{theorem}[subsection]{Theorem}

\newtheorem{proposition}[subsection]{Proposition}

\newtheorem{definition}[subsection]{Definition}

\newtheorem{claim}[subsection]{Claim}

\newtheorem{remark}[subsection]{Remark}

\def\loccitt{\emph{loc. cit.}}
\def\loccit{\emph{loc. cit. }}

\def\fsl{{\mathfrak{sl}}}
\def\fgl{{\mathfrak{gl}}}

\def\BA{{\mathbb{A}}}
\def\BC{{\mathbb{C}}}

\def\BN{{\mathbb{N}}}
\def\BF{{\mathbb{F}}}
\def\BP{{\mathbb{P}}}

\def\BQ{{\mathbb{Q}}}
\def\BZ{{\mathbb{Z}}}

\def\woo{\widehat{\otimes}}

\def\CA{{\mathcal{A}}}

\def\CE{{\mathcal{E}}}

\def\CS{{\mathcal{S}}}

\def\e{\varepsilon}

\def\and{\textrm{ }\&\textrm{ }}

\def\esym{\emph{Sym}}
\def\sym{\textrm{Sym}}

\def\sym{\textrm{Sym}}

\def\uu{{U_q(\dot{\fgl}_1)}}
\def\uup{{U_q^+(\dot{\fgl}_1)}}
\def\uum{{U_q^-(\dot{\fgl}_1)}}

\def\uul{{U_q^\leq(\dot{\fgl}_1)}}
\def\uug{{U_q^\geq(\dot{\fgl}_1)}}

\def\UU{{U_{q_1,q_2}(\ddot{\fgl}_1)}}

\def\UUl{{U^\leq_{q_1,q_2}(\ddot{\fgl}_1)}}
\def\UUg{{U^\geq_{q_1,q_2}(\ddot{\fgl}_1)}}

\def\la{{\lambda}}
\def\sq{{\square}}
\def\bsq{{\blacksquare}}

\def\lamu{{\lambda \backslash \mu}}
\def\lanu{{\lambda \backslash \nu}}
\def\munu{{\mu \backslash \nu}}

\def\bla{{\boldsymbol{\la}}}

\def\blamu{{\boldsymbol{\lamu}}}

\def\bn{\bar{n}}

\def\dR{\dot{R}}
\def\ddR{{\ddot{R}}}

\begin{document}

\title[The $R$-matrix of the quantum toroidal algebra]{\Large{\textbf{The $R$-matrix of the quantum toroidal algebra}}}

\author[Andrei Negu\cb t]{Andrei Negu\cb t}
\address{MIT, Department of Mathematics, Cambridge, MA, USA}
\address{Simion Stoilow Institute of Mathematics, Bucharest, Romania}
\email{andrei.negut@@gmail.com}

\maketitle

\begin{abstract} We consider the $R$-matrix of the quantum toroidal algebra of type $\fgl_1$, both abstractly and in Fock space representations. We provide a survey of a certain point of view on this object which involves the elliptic Hall and shuffle algebras, and show how to obtain certain explicit formulas.

\end{abstract}

\section{Introduction} 

\noindent The quantum toroidal algebra $\UU$ is quite a fascinating object: ubiquitous, but not completely belonging to a single area of mathematics and physics. It can be interpreted as the quantum affinization of the deformed Heisenberg algebra:
$$
\uu = U_q(\widehat{\fgl}_1)
$$
although since the latter is not of Drinfeld-Jimbo type, this interpretation is a bit ad-hoc. The quantum toroidal algebra was studied by Ding-Iohara (\cite{DI}) and Miki (\cite{M}), and appeared in numerous places in both the mathematical and physical literature, where it is sometimes known as the deformed $W_{1+\infty}$ algebra (see \cite{AKMMMOZ 1, AKMMMOZ 2, AKMMMOZ 3, AKMMSZ, B, BFMZZ, BMZ, FJMM, FHMZ, GG, Sm, Ts} and many other works). It is connected with geometric representation theory (\cite{FT, K-theory, SV}) and from there with the $q$-deformed Alday-Gaiotto-Tachikawa relations (\cite{AGT, AFHKSY, AY, MO, W, SV 3}). Last but not least, the quantum toroidal algebra is related to double affine Hecke algebras in type A (\cite{SV 2}). \\

\noindent The main purpose of this note is to exploit two other incarnations of $\UU$: the elliptic Hall algebra (\cite{BS, S}) and the double shuffle algebra (\cite{FHHSY, FO, Shuf}), in order to study the universal$^*$ $R$-matrix\footnote{The terminology ``universal$^*$ $R$-matrix" means that $\ddR$ differs from the actual universal $R$-matrix by certain powers of $q$, see \eqref{eqn:underlined} and \eqref{eqn:actual r matrix}, as shown in Section 2.2 of \cite{FJMM}. This is a well-known feature of quantum groups, where the analogous notion is the quasi $R$-matrix of \cite{L}.} :
$$
\ddR \in \UU \woo \UU 
$$ 
Using the tools developed in \cite{BS}, one obtains the following formula: \\

\begin{theorem}
\label{thm:main}

The universal$^*$ $R$-matrix can be factored as:
\begin{equation}
\label{eqn:factorization}
\ddR = \prod_{\text{coprime } (a,b) \in \BN \times \BZ \sqcup (0,1)} \exp \left[ \sum_{d=1}^\infty \frac {P_{da,db} \otimes P_{-da,-db}}d \frac {\left( q^{\frac d2} - q^{-\frac d2} \right)}{\left(q_1^{\frac d2} - q_1^{-\frac d2} \right) \left(q_2^{\frac d2} - q_2^{-\frac d2} \right)} \right] \qquad
\end{equation}
in terms of the generators $\{P_{n,m} \in \UU\}_{(n,m) \in \BZ^2 \backslash (0,0)}$ constructed by \cite{BS} and \cite{S} (see Subsection \ref{sub:eha} for a review). The product is taken in increasing order of $\frac ba$. \\

\end{theorem}

\noindent Formula \eqref{eqn:factorization} arises from the fact that products of the $P_{n,m}$'s in increasing order of slope form an orthogonal basis of the quantum toroidal algebra, which itself stems from the fact that coherent sheaves on an elliptic curve have Harder-Narasimhan filtrations. Moreover, \eqref{eqn:factorization} may be interpreted as a quantum toroidal version of the celebrated product formulas for universal $R$-matrices of quantum groups from \cite{D, KR, KT, LS, LSS, R}. The $\ddot{\fgl}_n$ analogue of the \eqref{eqn:factorization} was studied in \cite{Tor}. \\

\noindent Combining \eqref{eqn:factorization} with the shuffle algebra computations developed in \cite{Shuf}, one can obtain explicit formulas for the image of $\ddR$ in two types of Fock spaces\footnote{The representations $F^\uparrow_u$ and $F^\rightarrow_u$ were denoted by $F^{(0,1)}_u$ and $F^{(1,0)}_u$, respectively, in \cite{AFS}.}:
$$
\UU \curvearrowright F^\uparrow_u \text{ and }  F^\rightarrow_u 
$$
These formulas can be found in Theorems \ref{thm:right} and \ref{thm:up}, respectively, and can be used to understand tensor products of Fock spaces as representations of the quantum toroidal algebra. As explained in \cite{AFHKSY}, such tensor products govern the five-dimensional AGT relations (\cite{AGT, AY}). \\

\noindent I would like to thank Mikhail Bershtein, Jean-Emile Bourgine, Boris Feigin, Alexandr Garbali, Roman Gonin, Andrei Okounkov, Francesco Sala, Olivier Schiffmann, Junichi Shiarishi, Yan Soibelman and Alexander Tsymbaliuk for numerous wonderful discussions on the subject of quantum toroidal algebras over the years. I gratefully acknowledge NSF grants DMS-1760264 and DMS-1845034, as well as support from the Alfred P. Sloan Foundation. \\

\section{Warm-up: the deformed Heisenberg algebra}

\subsection{} Let $\BF$ be a field, implicitly the ground field of all our constructions. Recall that a bialgebra is an algebra $A$ with unit $1$ which is endowed with homomorphisms:
$$
\Delta : A \rightarrow A \otimes A \qquad \text{and} \qquad \e : A \rightarrow \BF 
$$
called \underline{coproduct} and \underline{counit}, respectively, which satisfy certain compatibility properties. We will often employ Sweedler notation for the coproduct:
\begin{equation}
\label{eqn:sweedler} 
\Delta(a) = a_1 \otimes a_2 
\end{equation}
the meaning of which is that there is an implied summation of tensors in the right-hand side. Given two bialgebras $A^+$ and $A^-$, a pairing between them:
\begin{equation}
\label{eqn:pair 0}
\langle \cdot, \cdot \rangle :  A^+ \otimes A^- \rightarrow \BF
\end{equation}
is called a \underline{bialgebra pairing} if it intertwines the product and coproduct as below:
\begin{align} 
&\langle a \cdot a', b \rangle = \left \langle a \otimes a', \Delta^{\text{op}}(b) \right \rangle \label{eqn:bialg 1} \\
&\langle a, b \cdot b' \rangle = \left \langle \Delta(a), b \otimes b' \right \rangle  \label{eqn:bialg 2}
\end{align}
for all $a,a' \in A^+$, and $b,b' \in A^-$. Given such a bialgebra pairing, we can form the \underline{Drinfeld double} (\cite{D}) of the bialgebras $A^+$ and $A^-$, namely the vector space:
\begin{equation}
\label{eqn:double}
A = A^+ \otimes A^-
\end{equation} 
One can make $A$ into a bialgebra by requiring that $A^+ \cong A^+ \otimes 1$ and $A^- \cong 1 \otimes A^-$ be sub-bialgebras, and that the multiplication of elements coming from different tensor factors be constrained by the relation:
\begin{equation}
\label{eqn:dd}
a_1 \cdot b_1 \langle a_2, b_2 \rangle = \langle a_1, b_1 \rangle b_2 \cdot a_2
\end{equation}
for all $a_1,b_1 \in A^+$ and $a_2,b_2 \in A^-$. \\

\begin{remark}

To define the bialgebra structure on \eqref{eqn:double} using relation \eqref{eqn:dd}, one needs all bialgebras involved to be \underline{Hopf} algebras. In other words, there must exist \underline{antipode} maps $S : A^\pm \rightarrow A^\pm$ which satisfy certain compatibility conditions with the product, coproduct and pairing, and then relation \eqref{eqn:dd} will be equivalent to:
\begin{equation}
\label{eqn:antipode}
\langle S^{-1}(a_1), b_1 \rangle a_2 b_2 \langle a_3, b_3 \rangle = b \cdot a 
\end{equation}
for all $a \in A^+$ and $b \in A^-$. Formula \eqref{eqn:antipode} is the one which allows to unambiguously define the product on the vector space \eqref{eqn:double}. The reason why we do not write down the antipode explicitly is that in all cases studied in the present paper, it exists and is uniquely determined by the bialgebra structure, and it is a straightforward exercise to write it down and to check that it satisfies all the required compatibility properties. \\

\end{remark}

\subsection{} 
\label{sub:r}

If the pairing \eqref{eqn:pair 0} is non-degenerate, then we may define:
\begin{equation}
\label{eqn:univ 0}
R = \sum_{i} a_i \otimes b_i \in A^+ \otimes A^- \hookrightarrow A \otimes A
\end{equation}
as $\{a_i, b_i\}_i$ go over any set of dual bases with respect to the pairing. The canonical tensor \eqref{eqn:univ 0} is called the \underline{universal $R$-matrix} of $A$, and satisfies the properties:
\begin{equation}
\label{eqn:univ prop 1}
R \cdot \Delta(a)  = \Delta^{\text{op}}(a) \cdot R
\end{equation}
for any $a \in A$ (here $\Delta^{\text{op}}$ denotes the opposite coproduct, which is obtained from $\Delta$ by switching the two tensor factors), as well as:
\begin{align}
&(\Delta \otimes 1)R = R_{13} R_{23} \label{eqn:univ prop 2} \\
&(1 \otimes \Delta)R = R_{13} R_{12} \label{eqn:univ prop 3}
\end{align}
where $R_{12} = R \otimes 1$, $R_{23} = 1 \otimes R$, and $R_{13}$ is defined analogously. The importance of this construction is the following: property \eqref{eqn:univ prop 1} implies that for any representations $V,W\in \text{Rep}(A)$, the operator $R^V_W$ given by: 
$$
A \otimes A \rightarrow \text{End}(V \otimes W), \quad R \mapsto R^V_{W}
$$ 
intertwines the $A$-module structures $V \otimes W$ and $W \otimes V$ (up to a swap of the factors). We may also perform this construction for a single representation $V \in \text{Rep}(A)$:
\begin{align*}
&A \otimes A \rightarrow A \otimes \text{End}(V), \qquad R \mapsto R_V \\
&A \otimes A \rightarrow \text{End}(V) \otimes A, \qquad R \mapsto R^V
\end{align*} 
Explicitly, the way one defines $R_V$ (respectively $R^V$) is to write $R$ as a sum of tensors $a\otimes b \in A \otimes A$, and replace each $b$ (respectively $a$) that appears in such tensors by the corresponding endomorphism of $V$ prescribed by the $A$-module structure of $V$. Given a vector and covector $v \in V$, $\lambda \in V^\vee$, we may therefore consider:
\begin{equation}
\label{eqn:notation}
_\lambda R_v \in A \qquad \text{and} \qquad ^\lambda R^v \in A
\end{equation}
obtained by taking the $\langle \lambda | v \rangle$ matrix coefficient of $R_V$ (respectively $R^V$) in the second (respectively first) tensor factor. \\

\subsection{} Consider two formal parameters $q_1$ and $q_2$, and set:
\begin{equation}
\label{eqn:q}
q = q_1 q_2 
\end{equation}
We will slightly abuse notation by writing $\BQ(q_1,q_2)$ instead of $\BQ(q_1^{\frac 12}, q_2^{\frac 12})$, which the reader should interpret as the fact that we fix square roots of $q_1$ and $q_2$. As these square roots are simply cosmetic, and not essential, features of the theory, this abuse seems acceptable. Let us consider the \underline{deformed Heisenberg} algebra:
$$
\uu = \BQ(q_1,q_2) \Big \langle P_n, c^{\pm 1} \Big \rangle_{n \in \BZ \backslash 0} \Big/^{c \text{ central}}_{\text{relation \eqref{eqn:heis}}}
$$
where:
\begin{equation}
\label{eqn:heis}
\Big[ P_n, P_{n'} \Big] = \delta_{n+n'}^0 \frac {n \left(q_1^{\frac n2} - q_1^{-\frac n2} \right) \left(q_2^{\frac n2} - q_2^{-\frac n2} \right)}{\left(q^{-\frac n2} - q^{\frac n2} \right)} \left(c^n - c^{-n} \right)
\end{equation}
It is a bialgebra with respect to the coproduct determined by:
$$
\Delta(c) = c \otimes c
$$
$$
\Delta(P_n) = \begin{cases} P_n \otimes 1 + c^n \otimes P_{n} &\text{if } n > 0 \\ P_n \otimes c^n + 1 \otimes P_n &\text{if } n < 0 \end{cases}
$$
and the counit determined by $\e(c) = 1$, $\e(P_n) = 0$ for all $n$. Moreover:
\begin{align*}
&\uug = \BQ(q_1,q_2)[P_n,c^{\pm 1}]_{n \in \BN} \\
&\uul = \BQ(q_1,q_2)[P_{-n},c^{\pm 1}]_{n \in \BN}
\end{align*} 
(here $\BN = \BZ_{>0}$) are sub-bialgebras of $\uu$, and there is a bialgebra pairing:
\begin{equation}
\label{eqn:pair heis}
\langle \cdot, \cdot \rangle : \uug \otimes \uul \rightarrow \BQ(q_1,q_2)
\end{equation} 
determined by:
\begin{equation}
\label{eqn:pair heis 1}
\langle c, - \rangle = \langle -, c\rangle = \e(-) 
\end{equation}
\begin{equation}
\label{eqn:pair heis 2}
\Big \langle P_n, P_{-n'} \Big \rangle = \delta_{n'}^n \frac {n \left(q_1^{\frac n2} - q_1^{-\frac n2} \right) \left(q_2^{\frac n2} - q_2^{-\frac n2} \right)}{\left( q^{\frac n2} - q^{-\frac n2} \right)}
\end{equation}
and properties \eqref{eqn:bialg 1}-\eqref{eqn:bialg 2}. It is easy to show that:
$$
\uu = \uug \otimes \uul \Big / (c \otimes 1 - 1 \otimes c)
$$ 
is the Drinfeld double constructed with respect to the pairing \eqref{eqn:pair heis}. \\

\begin{remark} Note that one can change the right-hand side of \eqref{eqn:heis} to any scalars that depend on $n$, and alternatively, this can be achieved by changing the right-hand side of \eqref{eqn:pair heis 2}. The reason why we prefer the scalars above is that they naturally appear in Macdonald polynomial theory (see \cite{Ops} for a brief survey of the connection) and in the study of the quantum toroidal algebra (see Section \ref{sec:tor}). \\
\end{remark}

\subsection{} It is easy to see that the restriction of the pairing \eqref{eqn:pair heis} to the subalgebras:
\begin{align*}
&\uup = \BQ(q_1,q_2)[P_n]_{n \in \BN} \\
&\uum = \BQ(q_1,q_2)[P_{-n}]_{n \in \BN}
\end{align*} 
is non-degenerate. Indeed, as $\bn = (n_1 \geq \dots \geq n_t)$ goes over partitions, the products:
$$
P_{\pm \bn} = P_{\pm n_1} \dots P_{\pm n_t}
$$
give rise to orthogonal bases with respect to \eqref{eqn:pair heis}:
$$
\Big \langle P_{\bn}, P_{-\bn'} \Big \rangle = \delta_{\bn'}^{\bn} z_{\bn} 
$$
where:
$$
z_{\bn} = \bn! \prod_{i=1}^t \frac {n_i \left(q_1^{\frac {n_i}2} - q_1^{-\frac {n_i}2} \right) \left(q_2^{\frac {n_i}2} - q_2^{-\frac {n_i}2} \right)}{\left( q^{\frac {n_i}2} - q^{-\frac {n_i}2} \right)}
$$
and $\bn!$ is the product of factorials of the number of times each positive integer $n$ appears in the partition $\bn$. One would like to invoke formula \eqref{eqn:univ 0} to conclude that the universal $R$-matrix of the deformed Heisenberg algebra is given by:
\begin{equation}
\label{eqn:r heis}
\dR := \sum_{\bn \text{ partition}}  \frac {P_{\bn} \otimes P_{-\bn}}{z_{\bn}} =  \exp \left[ \sum_{n=1}^\infty \frac {P_n \otimes P_{-n}}n \frac {\left( q^{\frac n2} - q^{-\frac n2} \right)}{\left(q_1^{\frac n2} - q_1^{-\frac n2} \right) \left(q_2^{\frac n2} - q_2^{-\frac n2} \right)} \right] \qquad 
\end{equation}
Because the exponential is an infinite sum, \eqref{eqn:r heis} lies in a certain completion:
$$
\dR \in \uu \woo \uu 
$$
We call \eqref{eqn:r heis} the \underline{universal$^*$ $R$-matrix} of $\uu$, and note that it slightly differs from the actual universal $R$-matrix because it is not equal to the canonical tensor of the pairing \eqref{eqn:pair heis}. The reason for this is that it does not account for the powers of the central element $c$. To make matters worse, the expressions $c^n - c^{n-1}$ lie in the kernel of the bialgebra pairing, making it degenerate, and thus not even allowing us to construct the canonical tensor. We will now show how to fix the issue. \\

\subsection{}
\label{sub:fix}

To construct the actual universal $R$-matrix, one needs to introduce an ``almost central" element $d$ into the deformed Heisenberg algebra, and consider:
\begin{equation}
\label{eqn:tilde heis}
\tilde{U}_q(\dot{\fgl}_1) = \uu \otimes_\BQ \BQ[d^{\pm 1}] \Big/ \Big( [c,d] = 0, \ d P_n = q^{-n} P_n d \Big)
\end{equation}
with coproduct $\Delta(d) = d \otimes d$. Moreover, we must replace \eqref{eqn:pair heis 1} by:
\begin{equation}
\label{eqn:new pair}
\langle c,d \rangle = \langle d,c \rangle = q, \qquad \langle c \text{ or } d ,P_n \rangle = \langle P_n, c \text{ or } d \rangle = 0
\end{equation}
for all $n \neq 0$. It is straightforward to check that this gives rise to a bialgebra pairing between the two halves of \eqref{eqn:tilde heis}. However, now we have a second problem, in that it is not clear to construct the canonical tensor on the infinite-dimensional vector space $\BQ(q_1,q_2)[c^{\pm 1}, d^{\pm 1}]$. To remedy this issue, we set:
$$
q = e^{\hbar}
$$
and work over $\BQ((\hbar))$ instead of over $\BQ(q)$. Then we replace the elements $c$ and $d$ by their logarithms $\gamma$ and $\delta$, explicitly defined by:
$$
c = e^{\hbar \gamma}, \quad d = e^{\hbar \delta}
$$
Set $\Delta(\gamma) = \gamma \otimes 1 + 1 \otimes \gamma$ and $\Delta(\delta) = \delta \otimes 1 + 1 \otimes \delta$, and define the pairing by:
$$
\langle \gamma, \delta \rangle = \langle \delta, \gamma \rangle = \frac 1{\hbar} 
$$
and all other pairings involving $\gamma$ and $\delta$ are set equal to 0. With this in mind, the canonical tensor restricted to the subalgebra generated by $\gamma, \delta$ takes the form:
$$
\sum_{n,n'=0}^\infty \gamma^n \delta^{n'} \otimes \delta^n \gamma^{n'} \cdot \frac {\hbar^{n+n'}}{n!n'!} = q^{\gamma \otimes \delta + \delta \otimes \gamma}
$$
Therefore, the correct formula for the universal $R$-matrix of $\uu$ is:
\begin{equation}
\label{eqn:underlined}
R_{\uu} = \dR \cdot q^{\log_q c \otimes \log_q d + \log_q d \otimes \log_q c}
\end{equation}
where $\dR \in \uu \woo \uu$ is defined in \eqref{eqn:r heis}. We stress once again the fact that in order to properly define the expression \eqref{eqn:underlined}, one needs to make all the modifications explained in the present Subsection: introduce the ``almost central" element $d$, work over power series in $\log q$ and replace the elements $c$ and $d$ by their logarithms in base $q$. Since the power of $q$ in \eqref{eqn:underlined} will always act by a simple operator in all representations we are concerned with, we will henceforth focus on providing formulas for $\dR$ (which will be the interesting part of the $R$-matrix for us). \\

\subsection{} The basic representation of $\uu$ is the \underline{Fock space}:
\begin{equation}
\label{eqn:fock}
F = \BQ(q_1,q_2)[p_1,p_2,\dots]
\end{equation}
with the action given by: 
$$
c \mapsto q^{\frac 12}, \qquad d \mapsto q^{\text{deg}}
$$
(here, $\deg$ denotes the grading on the polynomial ring which sets $\deg p_n = n$) and:
\begin{align}
&P_{-n} \mapsto \text{multiplication by } p_n \\
&P_{n} \ \ \mapsto - n  \left(q_1^{\frac n2} - q_1^{-\frac n2} \right) \left(q_2^{\frac n2} - q_2^{-\frac n2} \right) \frac {\partial}{\partial p_n}
\end{align}
The universal$^*$ $R$-matrix \eqref{eqn:r heis} in a tensor product of Fock modules is therefore:
\begin{equation}
\label{eqn:r heis fock}
\dR^{F}_{F} =  \sum_{\bn = (n_1 \geq \dots \geq n_t)} \frac {\partial}{\partial p_{\bn}} \otimes p_{\bn}  \cdot \frac {1}{\bn!} \prod_{i=1}^t \left( q^{-\frac {n_i}2} - q^{\frac {n_i}2} \right)
\end{equation}
where $p_{\bn} = p_{n_1} \dots p_{n_t}$ and $\frac {\partial}{\partial p_{\bn}} = \frac {\partial}{\partial p_{n_1}} \dots \frac {\partial}{\partial p_{n_t}}$. We have:
$$
\dR^{F}_{F} \in \text{End}(F \otimes F)
$$
because all but finitely many of the summands in \eqref{eqn:r heis fock} act trivially on any vector of $F \otimes F$, thus making formula \eqref{eqn:r heis fock} a well-defined endomorphism. This will be the case with all infinite sums that we will write in the present paper. \\






\section{The quantum toroidal algebra}
\label{sec:tor}

\subsection{} We will now consider the quantum toroidal algebra of type $\fgl_1$ (also known as the Ding-Iohara-Miki algebra). Consider the rational function:
\begin{equation}
\label{eqn:def zeta}
\zeta(x) = \frac {(1-xq_1)(1-xq_2)}{(1-x)(1-xq)}
\end{equation}
and the formal delta series $\delta(z) = \sum_{k \in \BZ} z^k$. \\

\begin{definition}

(\cite{DI, M}) Let:
$$
\UU = \BQ(q_1,q_2) \Big \langle e_k, f_k, h_{m}, c_1^{\pm 1}, c_2^{\pm 1} \Big \rangle_{k \in \BZ, m \in \BZ \backslash 0} \Big /^{c_1, c_2 \text{ central}}_{\text{relations \eqref{eqn:rel tor 1}-\eqref{eqn:rel tor 6}}}
$$
where we construct the power series $e(z) = \sum_{k \in \BZ} \frac {e_k}{z^k}$, $f(z) = \sum_{k \in \BZ} \frac {f_k}{z^k}$, and let:
\begin{equation}
\label{eqn:rel tor 1}
[h_m, h_{m'}] = \frac {\delta_{m+m'}^0  m \left(c_2^m - c_2^{-m} \right)}{\left(q_1^{\frac m2} - q_1^{-\frac m2} \right) \left(q_2^{\frac m2} - q_2^{-\frac m2} \right)\left(q^{-\frac m2} - q^{\frac m2} \right)} 
\end{equation}
\begin{align}
&[h_m, e_k] = e_{k+m} \cdot \begin{cases} 1 &\text{if } m>0 \\ -c_2^m &\text{if } m<0 \end{cases} \label{eqn:rel tor 2} \\
&[h_{m}, f_k] = f_{k+m} \cdot \begin{cases} 1 &\text{if } m<0 \\ -c_2^m &\text{if } m>0 \end{cases} \label{eqn:rel tor 3}
\end{align}
\begin{equation}
\label{eqn:rel tor 4}
e(z)e(w) \zeta \left(\frac zw \right) = e(w) e(z) \zeta \left(\frac wz \right)
\end{equation}
\begin{equation}
\label{eqn:rel tor 5}
f(z)f(w) \zeta \left(\frac wz \right) = f(w) f(z) \zeta \left(\frac zw \right)
\end{equation}
\begin{equation}
\label{eqn:rel tor 6}
[e_k, f_{k'}] = \frac {\left(q_1^{\frac 12} - q_1^{-\frac 12} \right) \left(q_2^{\frac 12} - q_2^{-\frac 12} \right)}{\left(q^{- \frac 12} - q^{ \frac 12} \right)} \left(\underbrace{\psi_{k+k'} c_1 c_2^{-k'}}_{\text{if }k+k' \geq 0} - \underbrace{\psi_{k+k'} c_1^{-1} c_2^{-k}}_{\text{if }k+k' \leq 0} \right)
\end{equation}
where the elements $\psi_m$ are defined by the generating series:
$$
\sum_{m=0}^\infty \psi_{\pm m} \cdot x^m = \exp \left[ \sum_{m=1}^\infty \frac {h_{\pm m}}m \cdot x^m \left(q_1^{\frac m2} - q_1^{-\frac m2}\right)\left(q_2^{\frac m2} - q_2^{-\frac m2}\right) \left(q^{\frac m2} - q^{- \frac m2} \right) \right] 
$$
To make sense of relations \eqref{eqn:rel tor 4} and \eqref{eqn:rel tor 5}, one clears denominators in the rational functions $\zeta$ and identifies the coefficients of $z^k w^l$ in the left and right-hand sides. \\

\end{definition}

\subsection{} We note that $\UU$ is a bialgebra, with coproduct:
\begin{equation}
\label{eqn:cop 0}
\Delta(c_1) = c_1 \otimes c_1 \qquad \Delta(c_2) = c_2 \otimes c_2 
\end{equation}
\begin{equation}
\label{eqn:cop 1}
\Delta(h_m) = \begin{cases} h_m \otimes 1 + c_2^m \otimes h_m &\text{if } m > 0 \\ h_m \otimes c_2^m + 1 \otimes h_m &\text{if } m < 0 \end{cases}
\end{equation} 
\begin{align}
&\Delta(e_k) = e_k \otimes 1 + \sum_{m=0}^\infty c_1 c_2^{k-m} \psi_{m} \otimes e_{k-m} \label{eqn:cop 2} \\
&\Delta(f_k) = 1 \otimes f_k + \sum_{m=0}^\infty f_{k+m} \otimes c_1^{-1} c_2^{k+m} \psi_{-m} \label{eqn:cop 3}
\end{align}
and counit determined by $\e(c_1) = \e(c_2) = 1$, $\e(e_k) = \e(f_k) = \e(h_m) = 0$. Note that the coproduct is defined in a topological sense, as it takes values in the completion:
\begin{equation}
\label{eqn:cop}
\Delta : \UU \rightarrow \UU \woo \UU
\end{equation}
It is easy to see that the coproduct preserves the subalgebras:
\begin{align*}
&\UUg = \BQ(q_1,q_2) \Big \langle e_k, h_m, c_1^{\pm 1}, c_2^{\pm 1} \Big \rangle_{k \in \BZ, m \in \BN} \\
&\UUl = \BQ(q_1,q_2) \Big \langle f_k, h_{-m}, c_1^{\pm 1}, c_2^{\pm 1} \Big \rangle_{k \in \BZ, m \in \BN}
\end{align*} 
of $\UU$, and it is well-known that we have a triangular decomposition:

\begin{equation}
\label{eqn:triangular toroidal}
\UU = \UUg \otimes \UUl \Big/(c_i \otimes 1 - 1 \otimes c_i)_{i \in \{1,2\}}
\end{equation}
Moreover, there exists a bialgebra pairing:
\begin{equation}
\label{eqn:pair tor}
\langle \cdot, \cdot \rangle : \UUg \otimes \UUl \rightarrow \BQ(q_1,q_2)
\end{equation}
determined by the assignments:
\begin{equation}
\label{eqn:pair exp 1}
\langle c_i, - \rangle = \langle -,  c_i \rangle = \e(-)
\end{equation} 
\begin{equation}
\label{eqn:pair exp 2}
\Big \langle e_k, f_{-k} \Big \rangle = \frac {\left(q_1^{\frac 12} - q_1^{-\frac 12} \right) \left(q_2^{\frac 12} - q_2^{-\frac 12} \right)}{\left( q^{\frac 12} - q^{-\frac 12} \right)}
\end{equation} 
\begin{equation}
\label{eqn:pair exp 3}
\Big \langle h_m, h_{-m} \Big \rangle = \frac {m}{\left(q_1^{\frac m2} - q_1^{-\frac m2} \right) \left(q_2^{\frac m2} - q_2^{-\frac m2} \right) \left( q^{\frac m2} - q^{-\frac m2} \right)} 
\end{equation}
(all other pairings between the generators $e_k, f_k, h_m$ are 0). Note that \eqref{eqn:triangular toroidal} is the Drinfeld double with respect to the datum above. Therefore, to construct and study the universal $R$-matrix of the quantum toroidal algebra, we must find dual bases with respect to the pairing \eqref{eqn:pair tor}. To achieve this, we now turn to another incarnation of the quantum toroidal algebra, namely the elliptic Hall algebra. \\

\subsection{} 
\label{sub:eha}

Let us consider the following half planes:
$$
\BZ_+^2 = \{(n,m) \in \BZ^2 \text{ s.t. } n>0 \text{ or } n=0,m>0\}
$$
$$
\BZ_-^2 = \{(n,m) \in \BZ^2 \text{ s.t. } n<0 \text{ or } n=0,m<0\}
$$

\begin{definition}
\label{def:eha}

(\cite{BS}) The elliptic Hall algebra is:
$$
\CA = \BQ(q_1,q_2) \Big \langle P_{n,m}, c_1^{\pm 1}, c_2^{\pm 1} \Big \rangle_{(n,m) \in \BZ^2 \backslash (0,0)} \Big /^{c_1, c_2 \text{ central}}_{\text{relations \eqref{eqn:relation 1}, \eqref{eqn:relation 2}}}
$$
where we impose the following relations:
\begin{equation}
\label{eqn:relation 1}
[P_{n,m}, P_{n',m'}] = \delta_{n+n'}^0 \frac {d \left(q_1^{\frac d2} - q_1^{-\frac d2} \right) \left(q_2^{\frac d2} - q_2^{-\frac d2} \right)}{\left(q^{-\frac d2} - q^{\frac d2} \right)} \left(c_1^n c_2^m - c_1^{-n} c_2^{-m} \right)
\end{equation}
if $nm'=n'm$ and $n>0$, with $d = \gcd(m,n)$. The second relation states that whenever $nm'>n'm$ and the triangle with vertices $(0,0), (n,m), (n+n',m+m')$ contains no lattice points inside nor on one of the edges, then we have the relation:
\begin{equation}
\label{eqn:relation 2}
[P_{n,m}, P_{n',m'}] = \frac {\left(q_1^{\frac d2} - q_1^{-\frac d2} \right) \left(q_2^{\frac d2} - q_2^{-\frac d2} \right)}{\left(q^{- \frac 12} - q^{ \frac 12} \right)} Q_{n+n',m+m'}
\end{equation}
$$
\cdot \ \begin{cases}
c_1^{-n'} c_2^{-m'} & \text{if } (n,m) \in \BZ_\pm^2, (n',m') \in \BZ_\mp^2, (n+n',m+m') \in \BZ_\pm^2 \\
c_1^n c_2^m & \text{if } (n,m) \in \BZ_\pm^2, (n',m') \in \BZ_\mp^2, (n+n',m+m') \in \BZ_\mp^2 \\
1 & \text{otherwise}
\end{cases}
$$
where $d = \gcd(n,m)\gcd(n',m')$ (by the assumption on the triangle, we note that at most one of the pairs $(n,m), (n',m'), (n+n',m+m')$ can fail to be coprime), and:
$$
\sum_{k=0}^{\infty} Q_{ka,kb} \cdot x^k = \exp \left[ \sum_{k=1}^\infty \frac {P_{ka,kb}}k \cdot x^k \left(q^{\frac k2} - q^{- \frac k2} \right) \right] 
$$
for all coprime integers $a,b$. Note that $Q_{0,0} = 1$. \\

\end{definition}

\begin{remark}

In the notation of \cite{BS}, we have: 
$$
P_{n,m} = \left(q_1^{\frac d2} - q_1^{-\frac d2} \right) \left(q_2^{\frac d2} - q_2^{-\frac d2} \right) u_{n,m}
$$
where $d = \gcd(n,m)$, as well as $c_1^n c_2^m = \kappa_{n,m}$. \\

\end{remark} 

\noindent As shown in \loccitt, when we specialize the parameters $q_1$ and $q_2$ to the Frobenius eigenvalues of an elliptic curve $\CE$ over the finite field $\BF_q$, and set $c_1,c_2 \mapsto 1$, the algebra $\CA$ matches a certain subalgebra of the Drinfeld double of the Hall algebra of the category of coherent sheaves over $\CE$. The fact that the group $SL_2(\BZ)$ acts on the derived category of coherent sheaves on $\CE$ translates into the fact that the universal cover of this group acts on the algebra $\CA$ by automorphisms (the reason one needs the universal cover is the presence of the central elements):
\begin{align}
&\widetilde{\gamma} \cdot P_{n,m} = P_{an+cm,bn+dm} \left(c_1^{an+cm}c_2^{bn+dm} \right)^\# \label{eqn:slz} \\
&\widetilde{\gamma} \cdot c_1 = c_1^a c_2^b, \qquad \widetilde{\gamma} \cdot c_2 = c_1^c c_2^d
\end{align} 
(see (6.16) of \cite{BS} for how to define the integer $\#$) where:
$$
\widetilde{\gamma} \in \widetilde{SL_2(\BZ)} \quad \text{is a lift of} \quad \gamma = \begin{pmatrix}
 a & c \\ b & d \end{pmatrix} \in SL_2(\BZ)
$$

\subsection{} It was shown in \cite{S} that there exists an algebra isomorphism:
\begin{equation}
\label{eqn:iso}
\UU \stackrel{\sim}\rightarrow \CA
\end{equation}
generated by:
$$
e_k \mapsto P_{1,k}, \quad f_k \mapsto P_{-1,k}, \quad h_m \mapsto \frac {P_{0,m}}{\left(q_1^{\frac m2} - q_1^{-\frac m2}\right)\left(q_2^{\frac m2} - q_2^{-\frac m2}\right)}
$$
This isomorphism allows us to transport the bialgebra structure from $\UU$ to $\CA$, as well as the decomposition \eqref{eqn:triangular toroidal} and the pairing \eqref{eqn:pair tor}:
\begin{equation}
\label{eqn:triangular elliptic}
\CA = \CA^\geq \otimes \CA^\leq \Big/(c_i \otimes 1 - 1 \otimes c_i)_{i \in \{1,2\}}
\end{equation}
\begin{equation}
\label{eqn:pair elliptic}
\langle \cdot, \cdot \rangle : \CA^\geq \otimes \CA^\leq \rightarrow \BQ(q_1,q_2)
\end{equation}
where we consider the subalgebras:
\begin{align} 
&\CA^\geq = \BQ(q_1,q_2) \Big \langle P_{n,m}, c_1^{\pm 1}, c_2^{\pm 1} \Big \rangle_{(n,m) \in \BZ_+^2} \label{eqn:geq} \\
&\CA^\leq = \BQ(q_1,q_2) \Big \langle P_{n,m}, c_1^{\pm 1}, c_2^{\pm 1} \Big \rangle_{(n,m) \in \BZ_-^2} \label{eqn:leq} 
\end{align}
of $\CA$. Moreover, in terms of the $P_{n,m}$ generators, the pairing takes the form:
\begin{equation}
\label{eqn:pair gens}
\Big \langle P_{n,m}, P_{n',m'} \Big \rangle = \delta_{n+n'}^0 \delta_{m+m'}^0 \frac {d \left(q_1^{\frac d2} - q_1^{-\frac d2} \right) \left(q_2^{\frac d2} - q_2^{-\frac d2} \right)}{\left( q^{\frac d2} - q^{-\frac d2} \right)}
\end{equation} 
where $d = \gcd(n,m)$, for all $m,m',n,n'$. Just like in the case of $\UU$, the decomposition \eqref{eqn:triangular elliptic} realizes $\CA$ as the Drinfeld double with respect to the data above. \\

\subsection{} The relevance of the generators $P_{n,m}$ is that ordered products of these elements give rise to an orthogonal basis of $\CA$. Explicitly, it was shown in \cite{BS} that:
\begin{equation}
\label{eqn:pbw}
P_{\pm v} = P_{\pm n_1, \pm m_1} \dots P_{\pm n_t, \pm m_t}
\end{equation}
for any \underline{convex} lattice path: 
\begin{equation}
\label{eqn:convex}
v = \Big\{(n_1,m_1), \dots ,(n_t,m_t)\Big\}, \quad \frac {m_1}{n_1} \leq \dots \leq \frac {m_t}{n_t}, \quad (n_i,m_i) \in \BZ^2_+
\end{equation}
(we identify convex paths up to permuting those edges $(n_i,m_i)$ with the same slope, which does not change the product \eqref{eqn:pbw} due to \eqref{eqn:relation 1}) give rise to a basis:
\begin{equation}
\label{eqn:basis of shuffle}
\CA^{\geq} = \left( \bigoplus_{v \text{ convex}} \BQ(q_1,q_2) \cdot P_{v} \right) \otimes_\BQ \BQ[c_1^{\pm 1}, c_2^{\pm 1}]
\end{equation}
as well as the analogous statement involving $\CA^\leq$ and $P_{-v}$. Moreover, these bases are orthogonal with respect to the pairing \eqref{eqn:pair elliptic} (see Proposition 5.7 of \cite{Shuf} for a proof, although it is already implicit from \cite{BS}):
\begin{equation}
\label{eqn:orthogonal convex}
\Big \langle P_{v}, P_{-v'} \Big \rangle = \delta_{v'}^v z_v 
\end{equation}
In the formula above, for a convex lattice path \eqref{eqn:convex} we define:
\begin{equation}
\label{eqn:zv}
z_v = v! \prod_{i=1}^t \frac {d_i \left(q_1^{\frac {d_i}2} - q_1^{-\frac {d_i}2} \right) \left(q_2^{\frac {d_i}2} - q_2^{-\frac {d_i}2} \right)}{\left( q^{\frac {d_i}2} - q^{-\frac {d_i}2} \right)}
\end{equation}
where we denote $d_i = \gcd(n_i, m_i)$, and $v!$ is the product of factorials of the number of times each vector $(n,m)$ appears in the path $v$. Therefore, by analogy with \eqref{eqn:r heis}, we call the following tensor the universal$^*$ $R$-matrix:
\begin{equation}
\label{eqn:r tor}
\ddR := \sum_{v \text{ convex}}  \frac {P_{v} \otimes P_{-v}}{z_v} = 
\end{equation}
$$
= \prod_{\text{coprime } (a,b) \in \BN \times \BZ \sqcup (0,1)} \exp \left[ \sum_{d=1}^\infty \frac {P_{da,db} \otimes P_{-da,-db}}d \cdot \frac {\left( q^{\frac d2} - q^{-\frac d2} \right)}{\left(q_1^{\frac d2} - q_1^{-\frac d2} \right) \left(q_2^{\frac d2} - q_2^{-\frac d2} \right)} \right] 
$$
where the product on the second line is taken in increasing order of $\frac ba$. Because of the isomorphism \eqref{eqn:iso}, we will refer to \eqref{eqn:r tor} as lying in either algebra:
$$
\ddR \in \CA \woo \CA \cong \UU \woo \UU
$$
As explained in Subsection \ref{sub:fix}, the actual universal $R$-matrix of $\UU$ is:
\begin{equation}
\label{eqn:actual r matrix}
R_{\UU} = \ddR \cdot q^{\sum_{i=1}^2 \log_q c_i \otimes \log_q d_i + \log_q d_i \otimes \log_q c_i}
\end{equation}
where $d_1, d_2$ are elements that one must add to the algebra $\UU$ (and work over $\BQ((\hbar_1, \hbar_2))$ instead of over $\BQ(q_1,q_2)$, where $\hbar_i = \log q_i$). We refer the reader to Section 2.2 of \cite{FJMM} for details as to the correct setup, and henceforth focus on \eqref{eqn:r tor}. \\

\subsection{} As is clear from \eqref{eqn:r tor}, understanding the generators $P_{n,m} \in \CA$ is key. To make them explicit, we turn to the shuffle algebra incarnation of $\CA$. Specifically, the following is a trigonometric degeneration of the $\ddot{\fgl}_1$ version of \cite{FO}. \\

\begin{definition}
\label{def:shuf} 

(\cite{FHHSY}) Consider the $\BQ(q_1,q_2)$-vector space:
\begin{equation}
\label{eqn:big}
V = \bigoplus_{k \geq 0} \BQ(q_1,q_2)(z_{1}, \dots ,z_{k})^{\esym}
\end{equation}
of rational functions which are symmetric in the variables $z_1, \dots ,z_k$, for any $k$. We endow $V$ with an algebra structure by the \underline{shuffle product}:
$$
R(z_{1}, \dots ,z_{k}) * R'(z_{1}, \dots ,z_{k'}) =
$$
\begin{equation}
\label{eqn:mult}
= \frac 1{k! k'!} \cdot \esym \left[R(z_{1}, \dots ,z_{k})R'(z_{k+1}, \dots ,z_{k+k'}) \prod_{i=1}^k \prod_{j = k+1}^{k+k'} \zeta \left( \frac {z_i}{z_j} \right) \right]
\end{equation}
where \esym \ denotes the symmetrization operator:
$$
\esym \left( R(z_1, \dots ,z_k) \right) = \sum_{\sigma \in S(k)} R(z_{\sigma(1)}, \dots ,z_{\sigma(k)})
$$
The \underline{shuffle algebra} $\CS \subset V$ is defined as the set of rational functions of the form:
\begin{equation}
\label{eqn:shuf}
R(z_{1}, \dots ,z_{k}) = \frac {r(z_{1}, \dots ,z_{k})}{\prod_{1 \leq i \neq j \leq k} (z_{i} - z_{j}q)} 
\end{equation}
where $r$ is a symmetric Laurent polynomial that satisfies the \underline{wheel conditions}:
\begin{equation}
\label{eqn:wheel}
r(z_1, \dots ,z_k) \Big|_{\left\{ \frac {z_1}{z_2}, \frac {z_2}{z_3}, \frac {z_3}{z_1} \right\} = \left\{q_1,q_2, \frac 1q \right\}} = 0
\end{equation}
Condition \eqref{eqn:wheel} is a trigonometric version of Condition 3 in $\S$1.3 of \cite{FO}. \\

\end{definition}

\subsection{} 
\label{sub:explicit shuffle} 

It was observed in \cite{FT, SV} that there are algebra homomorphisms:
\begin{align}
&\CA^- \stackrel{\Upsilon^-}\rightarrow \CS, \qquad \quad P_{-1,k} \mapsto z_1^k \label{eqn:shuf iso plus} \\
&\CA^+ \stackrel{\Upsilon^+}\rightarrow \CS^{\text{op}}, \qquad \ P_{1,k} \mapsto z_1^k \label{eqn:shuf iso minus}
\end{align}
where the subalgebras $\CA^\pm \subset \CA$ are defined by:
$$
\CA^\pm = \BQ(q_1,q_2) \Big \langle P_{n,m} \Big \rangle_{\pm n > 0, m \in \BZ} \\
$$
The maps \eqref{eqn:shuf iso plus} and \eqref{eqn:shuf iso minus} were shown in \cite{Shuf} to be isomorphisms, and the images of the generators $P_{n,m}$ under these maps were also computed:
\begin{equation}
\label{eqn:explicit pnm}
\Upsilon^\pm ( P_{\pm n,m} ) = q^{\frac {n-d}2} R_{n,m}(z_1, \dots ,z_n) 
\end{equation}
where for all $n \in \BN$ and $m\in \BZ$, we let $d = \gcd(n,m)$, $a = \frac nd$ and:
$$
R_{n,m} = \sym \left[ \frac {\prod_{i=1}^n z_i^{\left \lfloor \frac {im}n \right \rfloor - \left \lfloor \frac {(i-1)m}n \right \rfloor}}{\prod_{i=1}^{n-1} \left(1 - \frac {qz_{i+1}}{z_i} \right)} \sum_{s=0}^{d-1} q^{s} \frac {z_{a(d-1)+1} \dots z_{a(d-s)+1}}{{z_{a(d-1)} \dots z_{a(d-s)}}} \prod_{1 \leq i < j \leq n} \zeta \left( \frac {z_i}{z_j} \right) \right] 
$$
We have a triangular decomposition:
\begin{equation}
\label{eqn:triangular 1}
\CA = \CA^+ \otimes \CA^0 \otimes  \CA^-, \quad \CA^0 = \BQ(q_1,q_2)\Big \langle P_{0,m}, c_1^{\pm 1}, c_2^{\pm 1} \Big \rangle_{m \in \BZ \backslash 0}
\end{equation}
which matches the well-known triangular decomposition of $\UU$ under \eqref{eqn:iso}. \\

\section{Fock spaces}

\subsection{} Let us recall the bijection between partitions $\lambda = (\lambda_1 \geq \lambda_2 \geq \dots)$ and Young diagrams. The latter are sets of $1 \times 1$ boxes placed in the first quadrant of the plane, with $\lambda_1$ boxes placed on the first row, $\lambda_2$ boxes on the second row etc. For example, the following is the Young diagram associated to the partition $(4,3,1)$:

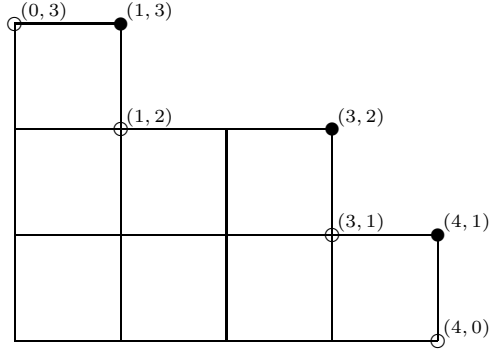
\begin{picture}(100,160)(-90,-15)
\label{fig}

\put(0,0){\line(1,0){160}}
\put(0,40){\line(1,0){160}}
\put(0,80){\line(1,0){120}}
\put(0,120){\line(1,0){40}}

\put(0,0){\line(0,1){120}}
\put(40,0){\line(0,1){120}}
\put(80,0){\line(0,1){80}}
\put(120,0){\line(0,1){80}}
\put(160,0){\line(0,1){40}}

\put(160,40){\circle*{5}}
\put(120,80){\circle*{5}}
\put(40,120){\circle*{5}}

\put(160,0){\circle{5}}
\put(120,40){\circle{5}}
\put(40,80){\circle{5}}
\put(0,120){\circle{5}}

\put(162,3){\scriptsize{$(4,0)$}}
\put(162,43){\scriptsize{$(4,1)$}}
\put(122,43){\scriptsize{$(3,1)$}}
\put(122,83){\scriptsize{$(3,2)$}}
\put(42,83){\scriptsize{$(1,2)$}}
\put(42,123){\scriptsize{$(1,3)$}}
\put(2,123){\scriptsize{$(0,3)$}}

\end{picture}

\noindent The hollow circles in the figure above will be called the \underline{inner corners} (abbreviated by ``i.c.") and the full circles will be called the \underline{outer corners} (abbreviated ``o.c.") of the partition. The \underline{weight} of a box is defined as the quantity:
\begin{equation}
\label{eqn:weight}
\chi_\sq = u q_1^x q_2^y
\end{equation}
where $(x,y)$ are the coordinates of the bottom left corner of the box, and $u$ is a parameter. Given two Young diagrams, we will write $\mu \subset \lambda$ if $\mu$ is contained in $\lambda$. If this happens, and $R$ is a symmetric rational function in $|\lamu|$ variables, write:
\begin{equation}
\label{eqn:specialize}
R(\lamu) = R(\dots ,\chi_\sq,\dots)_{\sq \in \lambda \backslash \mu} 
\end{equation}
We set $R(\lamu) = 0$ if $\mu \not \subset \lambda$ or if $|\lamu|$ is not equal to the number of variables of $R$. \\

\begin{definition}
\label{def:action right}

(\cite{FT, SV}, see also \cite{K-theory}) Let $F_u^\uparrow$ be a vector space with a basis $|\lambda\rangle$ indexed by partitions. Then the following formulas determine an action $\CA \curvearrowright F_u^\uparrow$:
\begin{equation}
\label{eqn:action right 1}
c_1 \mapsto q^{\frac 12}, \qquad c_2 \mapsto 1,
\end{equation}
\begin{equation}
\label{eqn:action right 2}
\langle \mu | P_{0,\pm m} | \lambda \rangle = \pm \delta_\lambda^\mu q^{\mp \frac m2} \left( \sum_{\sq \text{ i.c. of }\lambda} \chi_\sq^{\pm m} - \sum_{\sq \text{ o.c. of }\lambda} \chi_\sq^{\pm m} \right)
\end{equation}
\footnote{Given an endomorphism $S$ of $F_u^\uparrow$, we write $\langle \mu| S |\lambda \rangle$ for the coefficient of $|\mu \rangle$ in $S \left(|\lambda \rangle\right)$.} and for all $X \in \CA^-$ with $\Upsilon^-(X) = R$ (resp. $Y \in \CA^+$ with $\Upsilon^+(Y) = R$), we have:
\begin{align} 
&\langle \lambda | X | \mu \rangle = R(\lamu) \cdot \sigma^{|\lamu|} \prod_{\bsq \in \lamu} \frac {\prod_{\sq \text{ o.c. of }\lambda} \left(1 - \frac {\chi_\sq}{\chi_\bsq} \right)}{\prod_{\sq \text{ i.c. of }\lambda} \left(1 - \frac {\chi_\sq}{\chi_\bsq} \right)}  \label{eqn:action right 3} \\ 
&\langle \mu | Y | \lambda \rangle = R(\lamu) \cdot \bar{\sigma}^{|\lamu|} \prod_{\bsq \in \lamu} \frac {\prod_{\sq \text{ i.c. of }\mu} \left(1 - \frac {\chi_\sq}{q\chi_\bsq} \right)}{\prod_{\sq \text{ o.c. of }\mu} \left(1 - \frac {\chi_\sq}{q\chi_\bsq} \right)}  \label{eqn:action right 4} 
\end{align}
where $\sigma = \frac {(1-q_1)(1-q_2)}{1-q}$ and $\bar{\sigma} = \sigma q^{\frac 12}$. \\

\end{definition}

\noindent The reader might ask how to interpret the evaluation $R(\lamu)$, given that elements $R$ of the shuffle algebra take the form \eqref{eqn:shuf}, and thus have poles at $z_i - z_j q$. The answer lies in the wheel conditions. One first defines the specialization:
$$
\rho(y_1,y_2,\dots) = R \left(y_1 q_1^{\lambda_1 - \mu_1}, \dots ,y_1 q_1^{\lambda_1-1}, y_2 q_1^{\lambda_2 - \mu_2}, \dots ,y_2 q_1^{\lambda_2-1}, \dots \right)
$$
(which is allowed, because $R$ has no poles at $z_i - z_j q_1$) and then invoke the wheel conditions \eqref{eqn:wheel} to conclude that $\rho$ has no poles at $y_i - y_j q_2$. Then we define:
$$
R(\lamu) = \rho(1,q_2,q_2^2,\dots)
$$

\begin{remark} 

If we replace individual partitions $\lambda$ by $r$-tuples of partitions $\bla = (\lambda^1,\dots,\lambda^r)$, then straightforward analogues of formulas \eqref{eqn:action right 1}-\eqref{eqn:action right 4} yield an action:
$$
\CA \curvearrowright F_{u_1}^\uparrow \otimes \dots \otimes F_{u_r}^\uparrow
$$
To this end, one must replace the weight \eqref{eqn:weight} of a box in an individual partition by the weight $\chi_\sq = u_i q_1^x q_2^y$ of a box $\sq$ located at coordinates $(x,y)$ in the $i$-th constituent partition of an $r$-tuple of partitions $\bla$. We refer the reader to \cite{K-theory} for details, and for the connection to moduli spaces of rank $r$ sheaves on the affine plane. \\

\end{remark} 
	
\subsection{} 
\label{sub:right}

Let us use the notation $F_u^\rightarrow = \BQ(q_1,q_2)[p_1,p_2,\dots]$ for the Fock space \eqref{eqn:fock}. \\

\begin{definition}
\label{eqn:action up} 

(\cite{FHHSY, FJMM}) The following formulas determine an action $\CA \curvearrowright F_u^\rightarrow$:
\begin{equation}
\label{eqn:action up 1}
c_1 \mapsto 1, \qquad c_2 \mapsto q^{\frac 12},
\end{equation}
\begin{align}
&P_{0,-m} \mapsto \text{multiplication by } p_m  \label{eqn:action up 2} \\
&P_{0,m} \ \ \mapsto - m \left(q_1^{\frac m2} - q_1^{-\frac m2} \right) \left(q_2^{\frac m2} - q_2^{-\frac m2} \right) \frac {\partial}{\partial p_m} \label{eqn:action up 3}
\end{align}
and $\forall X \in \CA^-$ with $\Upsilon^-(X) = R(z_1,\dots,z_n)$ (respectively $Y \in \CA^+$ with $\Upsilon^+(Y) = R$):
\begin{equation}
\label{eqn:action up 4}  
X \mapsto \frac {(u q^{-\frac 12})^n}{n!} \int^{|q_1|,|q_2| < 1}_{|z_1| = \dots = |z_n| = 1} \frac {R(z_1,\dots,z_n)}{\prod_{1\leq i \neq j \leq n} \zeta \left( \frac {z_i}{z_j} \right)}  \prod_{a=1}^n  \frac {dz_a}{2\pi i z_a}
\end{equation} 
$$
\exp \left[\sum_{k=1}^\infty \frac {z_1^k+\dots+z_n^k}k \cdot q^{\frac k2} p_k \right] \exp \left[-\sum_{k=1}^\infty (z_1^{-k}+\dots+z_n^{-k}) \cdot q^{-\frac k2} (1-q_1^k)(1-q_2^k)  \frac {\partial}{\partial p_k} \right] 
$$
\begin{equation}
\label{eqn:action up 5}  
Y \mapsto \frac {(-u^{-1}q^{\frac 12})^{n}}{n!} \int^{|q_1|,|q_2| > 1}_{|z_1| = \dots = |z_n| = 1} \frac {R(z_1,\dots,z_n)}{\prod_{1\leq i \neq j \leq n} \zeta \left( \frac {z_i}{z_j} \right)} \prod_{a=1}^n  \frac {dz_a}{2\pi i z_a}  
\end{equation} 
$$
\exp \left[-\sum_{k=1}^\infty \frac {z_1^k+\dots+z_n^k}k \cdot p_k \right] \exp \left[\sum_{k=1}^\infty (z_1^{-k}+\dots+z_n^{-k}) (1-q_1^{-k})(1-q_2^{-k}) \cdot \frac {\partial}{\partial p_k} \right] 
$$
The integrals in \eqref{eqn:action up 4} and \eqref{eqn:action up 5} are contour integrals, and the parameters $q_1$ and $q_2$ are interpreted as complex numbers satisfying the conditions displayed in the superscripts of the integral signs. 

\end{definition}

\noindent We have an isomorphism of vector spaces:
\begin{equation}
\label{eqn:iso fock}
\Psi : F_u^\uparrow \stackrel{\sim}\longrightarrow F_u^\rightarrow
\end{equation}
obtained by sending $|\lambda\rangle$ to the modified Macdonald polynomial associated to the partition $\lambda$ (see \cite{GH, H}) for parameters $(q,t) \leftrightarrow (q_1^{-1},q_2^{-1})$. However, the isomorphism \eqref{eqn:iso fock} also respects the $\CA$ actions, up to rotation by $90$ degrees:
\begin{equation}
\label{eqn:iso actions}
\Psi(P_{n,m} \cdot x) = q^{\frac {m \varepsilon_{n,m}}2}P_{-m,n} \cdot \Psi(x)
\end{equation} 
for any $x \in F_u^\uparrow$ and $(n,m) \in \BZ^2 \backslash (0,0)$, where $\varepsilon_{n,m}$ is defined to be $-1$ if $(n,m)$ lies in the second or fourth quadrant (including the horizontal axis, but excluding the vertical axis) and 0 otherwise. To prove \eqref{eqn:iso actions}, note that the algebra $\CA$ is generated by $P_{n,m}$ with $(n,m) \in (\pm 1, 0), (0,\pm 1)$. Moreover:
$$
P_{n,m} \mapsto (c_1^n c_2^m)^{\e_{n,m}} P_{-m,n}
$$
is an algebra automorphism, namely the particular case of \eqref{eqn:slz} for:
\begin{equation}
\label{eqn:90}
\gamma = \begin{pmatrix}
0 & -1 \\ 1 & 0 \end{pmatrix} \in SL_2(\BZ)
\end{equation}
Therefore, in order to prove \eqref{eqn:iso actions} for a general lattice point $(n,m)$, it suffices to prove it for the four special lattice points $(\pm 1, 0), (0,\pm 1)$. All four of these statements are well-known facts in Macdonald polynomial theory. \\

\begin{remark}
\label{rem:rg}
The gist of \eqref{eqn:iso fock} and \eqref{eqn:iso actions} is that $F_u^\uparrow$ and $F_u^\rightarrow$ can be perceived as the same module, up to the automorphism provided by the matrix \eqref{eqn:90}. One could turn the problem on its head, and define for any $\gamma \in SL_2(\BZ)$ the module:
\begin{equation}
\label{eqn:gamma}
F_u^\gamma = F_u^\rightarrow
\end{equation}
but any element $a \in \CA$ acts on the left-hand side just like $\gamma(a) \in \CA$ acts on the right-hand side, where $\gamma(a)$ denotes the automorphism \eqref{eqn:slz} for some lift of $\gamma$. Up to a simple isomorphism (namely conjugation by powers of the famous $\nabla$ operator, see \cite{BG}), this module structure only depends on:
$$
\gamma \begin{pmatrix} 1 \\ 0 \end{pmatrix} = \begin{pmatrix} a \\ b \end{pmatrix}  \in \BZ^2
$$
i.e. the choice of $b/a \in \BQ \sqcup \infty$. It is an interesting open problem to present the module structure on $F_u^\gamma$ in such a way that any $X \in \CA^\pm$ acts by an explicit formula that essentially involves only the rational function $R = \Upsilon^\pm(X)$. This was provided in \eqref{eqn:action right 3}-\eqref{eqn:action right 4} for the vertical slope and in \eqref{eqn:action up 4}-\eqref{eqn:action up 5} for the horizontal slope. Such rotated modules for the quantum toroidal algebra play an important role in the study of refined link invariants, \cite{GN 1, GN 2}. \\

\end{remark}
	
\subsection{} 
\label{sub:up action}

We will now compute the matrix coefficients of the universal$^*$ $R$-matrix \eqref{eqn:r tor} in the two types of Fock spaces considered in the present paper. Recall that convex paths go over $\BZ_+^2$, and we will use the term \underline{convex}$^*$ path to refer to those whose edges do not point directly up (i.e. we restrict to $n_i > 0$ in \eqref{eqn:convex}). The \underline{size} of such a path is the $x$-coordinate of the lattice point where it terminates, i.e. the number $n_1+\dots+n_t$ in \eqref{eqn:convex}. Any convex path $v$ can be obtained by concatenating a convex$^*$ path $v^*$ with a vertical path, hence we can uniquely write:
$$
P_v = P_{v^*} P_{0,\bar{n}}
$$
where $P_{0,\bar{n}} = P_{0,n_1} \dots P_{0,n_t}$ for any partition $\bar{n} = (n_1 \geq \dots \geq n_t)$. Therefore:
\begin{equation}
\label{eqn:new r tor 1}
\boxed{\ddR = R' R''} 
\end{equation}
where:
\begin{align}
&R' = \sum_{v \text{ convex}^*}  \frac {P_{v} \otimes P_{-v}}{z_v} \in \CA^+ \woo \CA^- \cong \CS^{\text{op}} \woo \CS \label{eqn:new r tor 2} \\
&R'' = \sum_{\bn \text{ partition}}  \frac {P_{0,\bn} \otimes P_{0,-\bn}}{z_{\bn}} \in \CA^0 \woo \CA^0 \label{eqn:new r tor 3}
\end{align} 
We will now compute the matrix coefficients of $R'$ and $R''$ in the module $F^\uparrow_u$. Take any two partitions $\mu$ and $\lambda$, and regard them as a vector and covector:
$$
\langle \lambda | \in (F_u^\uparrow)^\vee \qquad |\mu \rangle \in F_u^\uparrow
$$
Then we have (recall \eqref{eqn:notation} for notation pertaining to matrix coefficients of $R', R''$ in one of the two tensor factors):
\begin{align} 
&_{\langle \lambda |}R'_{|\mu \rangle} =  \sum_{v \text{ convex}^*}  \frac {P_v}{z_v} \cdot \langle \lambda | P_{-v} |\mu \rangle \label{eqn:prime} \\
&_{\langle \lambda |}R''_{|\mu \rangle} =  \sum_{\bn \text{ partition}}  \frac {P_{0,\bn}}{z_{\bn}} \langle \lambda | P_{0,-\bn} |\mu \rangle \label{eqn:double prime}
\end{align}
By formulas \eqref{eqn:action right 2} and \eqref{eqn:action right 3}, the right-hand sides above are non-zero only if $\mu \subset \lambda$. \\

\begin{claim} 
\label{claim:first}
	
The right-hand side of \eqref{eqn:double prime} is equal to:
\begin{equation}
\label{eqn:claim first}
\delta_\lambda^\mu \cdot \exp \left[ \sum_{k=1}^\infty \frac {P_{0,k}}k \cdot q^{\frac k2} \left( \sum_{\sq \text{ o.c. of }\mu} \chi_\sq^{-k} - \sum_{\sq \text{ i.c. of }\mu} \chi_\sq^{-k} \right) \right] 
\end{equation}

\end{claim} 

\noindent The claim above is a simple exercise, which we leave to the interested reader. \\

\noindent Meanwhile, \eqref{eqn:action right 3} implies that the right-hand side of \eqref{eqn:prime} is equal to:
\begin{equation}
\label{eqn:first}
\left[ \sum_{v \text{ convex}^* \text{ of size }|\lamu|}  \frac {P_v}{z_v} \cdot  \Upsilon^-(P_{-v})(\lamu) \right] \sigma^{|\lamu|} \prod_{\bsq \in \lamu} \frac {\prod_{\sq \text{ o.c. of }\lambda} \left(1 - \frac {\chi_\sq}{\chi_\bsq} \right)}{\prod_{\sq \text{ i.c. of }\lambda} \left(1 - \frac {\chi_\sq}{\chi_\bsq} \right)}  \qquad  \quad
\end{equation}
We note that expression \eqref{eqn:first} is an infinite sum of $P_v$'s, over convex$^*$ paths starting at 0 and ending at a point on the line $x = |\lamu|$. Therefore, it can only act on graded modules where such $P_v$'s act locally nilpotently, which will henceforth be called \underline{good} modules (for example, $F_{u'}^\rightarrow$ is good). Then let us consider the expression:
$$
W(x) = \sum_{k \in \BZ} P_{1,k}x^k
$$
In any good module, the expression $W(x_1) \dots W(x_n)$ makes sense when expanded in $|x_1| \gg \dots \gg |x_n|$. Moreover, it is clear from the shuffle algebra incarnation that:
\begin{equation}
\label{eqn:normal}
W(x_1,\dots, x_n) = W(x_1) \dots W(x_n) \prod_{1 \leq i < j \leq n} \zeta \left(\frac {x_j}{x_i} \right)
\end{equation}
is a symmetric expression in $x_1,\dots,x_n$, at least formally. In good representations, this means that the expression \eqref{eqn:normal} has the property that all its matrix coefficients are symmetric rational functions in $x_1,\dots,x_n$. \\

\begin{remark}

In the good representation $F_{u'}^\rightarrow$, the expression $W(x_1,\dots,x_n)$ acts by:
$$
\exp \left[-\sum_{k=1}^\infty \frac {x_1^{-k}+\dots+x_n^{-k}}k p_k \right] \exp \left[\sum_{k=1}^\infty (x_1^{k}+\dots+x_n^{k}) (1-q_1^{-k})(1-q_2^{-k})\frac {\partial}{\partial p_k} \right]
$$
times $(-{u'}^{-1} q^{\frac 12})^n$. \\

\end{remark}

\noindent As a consequence of formula (7.7) of \cite{W}, we have the following formula for the pairing:
\begin{equation}
\label{eqn:pairing shuffle 1}
\Big \langle W(x_1,\dots,x_n), a \Big \rangle = \Upsilon^-(a)(x_1,\dots,x_n) \qquad \forall a \in \CA^-
\end{equation}
Letting $a = P_{-v}$ for any convex$^*$ path $v$, and recalling that such convex paths give rise to orthogonal bases \eqref{eqn:orthogonal convex}, we obtain the following formula:
\begin{equation}
\label{eqn:formula}
W(x_1, \dots, x_n) = \sum_{v \text{ convex}^* \text{ of size }n} \frac {P_v}{z_v} \cdot \Upsilon^-(P_{-v})(x_1,\dots,x_n)
\end{equation}
If we plug this formula in \eqref{eqn:first}, then we conclude: \\

\begin{theorem}
\label{thm:right}

For any partitions $\lambda, \mu$, we have the following identity in $\CA^+$:
\begin{equation}
\label{eqn:first 3}
_{\langle \lambda |}R'_{|\mu \rangle} = \underline{W(\dots,\chi_\bsq,\dots)_{\bsq \in \lamu}} \cdot \sigma^{|\lamu|} \prod_{\bsq \in \lamu} \frac {\prod_{\sq \text{ o.c. of }\lambda} \left(1 - \frac {\chi_\sq}{\chi_\bsq} \right)}{\prod_{\sq \text{ i.c. of }\lambda} \left(1 - \frac {\chi_\sq}{\chi_\bsq} \right)} 
\end{equation}
Meanwhile, $_{\langle \lambda |}R''_{|\mu \rangle}$ is given by \eqref{eqn:claim first}. \\

\end{theorem}

\noindent In the module $F_{u'}^\rightarrow$, the underlined term in \eqref{eqn:first 3} is precisely the normal ordered product (3.12) of \cite{AFS}, which plays a key role in the construction of the intertwiner:
$$
F^\uparrow_u \otimes F^\rightarrow_{u'} \longrightarrow F^{\nearrow}_{-uu'}
$$
where $F^{\nearrow}_{-uu'}$ is the module \eqref{eqn:gamma} for:
$$
\gamma = \begin{pmatrix} 1 & 0 \\ 1 & 1 \end{pmatrix}
$$
This module was denoted by $F^{(1,1)}_{-uu'}$ in \loccitt \\

\begin{remark}

A similar formula to \eqref{eqn:first 3} holds if we replace the module $F_u^\uparrow$ by the MacMahon module of \cite{FJMM 2}, in which case the variables of the underlined term should go over the set of weights of boxes of a skew 3-dimensional partition $\blamu$. The generalization is straightforward, but the product of factors in \eqref{eqn:first 3} is more involved in the MacMahon case, and so we leave the details as an exercise. \\

\end{remark}

\subsection{}
\label{sub:right r-matrix}  

We will now consider the Fock space of Subsection \ref{sub:right}, and identify it with:
$$
F_u^\rightarrow = \BQ(q_1,q_2)[p_1,p_2,\dots] \stackrel{\sim}{\longrightarrow}\BQ(q_1,q_2)[x_1,x_2,\dots]^{\sym} = \Lambda
$$
via $p_n = x_1^n+x_2^n+\dots$. We note that a linear basis of $F_u^\rightarrow \cong \Lambda$ is given by:
\begin{equation}
\label{eqn:basis}
p_{\bn} = p_{n_1} \dots p_{n_t} 
\end{equation}
as $\bn = (n_1 \geq \dots \geq n_t)$ goes over partitions. We regard elements of $F_u^\rightarrow \cong \Lambda$ as functions $f[X]$, where $X$ is shorthand for the variable set $x_1,x_2,\dots$. We adopt ``plethystic notation", according to which one defines, for any symbol $z$:
\begin{equation}
\label{eqn:plethysm}
f[X \pm z] \in \Lambda [[z^{\pm 1}]]
\end{equation}
to be the image of $f[X]$ under the ring homomorphism $\Lambda \rightarrow \Lambda[[z^{\pm 1}]]$ that sends:
\begin{equation}
\label{eqn:plethysm 2}
p_n \mapsto p_n \pm z^n
\end{equation}
In other words, the way one computes \eqref{eqn:plethysm} is to expand $f[X]$ in the basis \eqref{eqn:basis}, and then replace each $p_n$ therein according to \eqref{eqn:plethysm 2}. Morally, the plethysm \eqref{eqn:plethysm} means ``add/remove $z$ from the list of variables", hence the notation. \\

\begin{claim}
\label{claim:plethysm}

For any $f[X] \in \Lambda$ and any variables $z_1,\dots,z_n$, we have:
\begin{multline}
\exp \left[\sum_{k=1}^\infty (z_1^{-k}+\dots+z_n^{-k}) (1-q_1^{-k})(1-q_2^{-k})\frac {\partial}{\partial p_k} \right] \cdot f[X] = \\ = f \left[ X + \left(1-q_1^{-1} \right)\left(1- q_2^{-1} \right) \sum_{i=1}^n z_i^{-1} \right] \label{eqn:pleth exp}
\end{multline}

\end{claim}

\noindent It is enough to prove the claim above for:
$$
f[X] = \prod_{i=1}^{\infty} \prod_{j=1}^t (1-x_i a_j) = \exp \left[- \sum_{k=1}^\infty \frac {a_1^k + \dots + a_t^k}k \cdot p_k \right]
$$
because the coefficients of such expressions in the variables $a_i$ provide a linear basis of $\Lambda$. The computation of \eqref{eqn:pleth exp} for $f$ as in the formula above is a straightforward exercise, which we leave to the interested reader. The following is also obvious. \\

\begin{claim}
\label{claim:adjoint}

The operators $-\frac {p_n}n$ and $(1-q_1^{-n})(1-q_2^{-n}) \frac {\partial}{\partial_{p_n}}$ are adjoint with respect to the (modified) Macdonald inner product:
\begin{equation}
\label{eqn:macdonald}
\langle \cdot, \cdot \rangle : \Lambda \otimes \Lambda \rightarrow \BQ(q_1,q_2)
\end{equation}
$$
\Big \langle p_{\bn}, p_{\bn'} \Big \rangle = \delta_{\bn'}^{\bn} \bn! \prod_{i=1}^t \Big[ - n_i (1-q_1^{-n_i})(1-q_2^{-n_i}) \Big]
$$

\end{claim}

\noindent We will now use the language above to compute the universal$^*$ $R$-matrix in the module $F_u^\rightarrow$. Let us recall the decomposition \eqref{eqn:new r tor 1}. Then we have:
$$
^{\langle f|}{R''}^{|g \rangle} = \sum_{\bn \text{ partition}}  \frac {P_{0,-\bn}}{z_{\bn}} \cdot \langle f , P_{0,\bn} g \rangle
$$
According to Claim \ref{claim:adjoint} and formula \eqref{eqn:action up 3}, we obtain (recall \eqref{eqn:notation} for notation pertaining to matrix coefficients of $R''$ in one of the two tensor factors):
\begin{equation}
\label{eqn:second}
^{\langle f|}{R''}^{|g \rangle} = \sum_{\bn \text{ partition}}  \frac {P_{0,-\bn}}{z_{\bn}} \cdot q^{\frac {|\bn|}2} \langle p_{\bn} f , g \rangle = \left \langle \exp \left[ \sum_{n=1}^\infty \frac {P_{0,-n}q^{\frac n2}}{n} \cdot p_n \right] \cdot f, g \right \rangle \qquad
\end{equation}
where $p_n$ acts on symmetric functions as multiplication by $p_n$, while the symbols $P_{0,-n} \in \CA^0$ are unaffected by their interaction with $f$ and $g$, or the pairing. \\

\noindent To compute the matrix coefficients of $R'$, we need to understand the action:
$$
\CA^+ \cong \CS^{\text{op}} \curvearrowright F_u^\rightarrow \cong \Lambda
$$ 
of \eqref{eqn:action up 5} in the language of the present Subsection. Claims \ref{claim:plethysm} and \ref{claim:adjoint} imply:
\begin{equation}
\label{eqn:new action}
\langle f | Y | g \rangle = \frac {(-u^{-1}q^{\frac 12})^{n}}{n!} \int^{|q_1|,|q_2| > 1}_{|z_1| = \dots = |z_n| = 1} \frac {\Upsilon^+(Y)(z_1,\dots,z_n)}{\prod_{1\leq i \neq j \leq n} \zeta \left( \frac {z_i}{z_j} \right)}  \prod_{a=1}^n  \frac {dz_a}{2\pi i z_a}
\end{equation}
$$
\left \langle f \left[ X + \left(1-q_1^{-1} \right)\left(1- q_2^{-1} \right) \sum_{i=1}^n z_i \right], g \left[ X + \left(1-q_1^{-1} \right)\left(1- q_2^{-1} \right) \sum_{i=1}^n z_i^{-1} \right] \right \rangle
$$
for any $Y \in \CA^+$, $f = f[X], g = g[X] \in \Lambda$, where the inner product is given by \eqref{eqn:macdonald}. Let us explain the gist of \eqref{eqn:new action}: for any two fixed symmetric polynomials $f$ and $g$, the second line of \eqref{eqn:new action} is a Laurent polynomial in the variables $z_1,\dots,z_n$, which one then integrates against the rational function on the first line of \eqref{eqn:new action}. \\

\noindent With this in mind, we may use the definition of $R'$ in \eqref{eqn:new r tor 2} to compute:
$$
^{\langle f|}{R'}^{|g\rangle} = \sum_{n=0}^\infty  \frac {(-u^{-1}q^{\frac 12})^{n}}{n!} \sum_{v \text{ convex}^* \text{ of size }n}  \frac {P_{-v}}{z_v} \int^{|q_1|,|q_2| > 1}_{|z_1| = \dots = |z_n| = 1} \frac {\Upsilon^+(P_v)(z_1,\dots,z_n)}{\prod_{1\leq i \neq j \leq n} \zeta \left( \frac {z_i}{z_j} \right)}  \prod_{a=1}^n  \frac {dz_a}{2\pi i z_a}
$$
$$
\left \langle f \left[ X + \left(1-q_1^{-1} \right)\left(1- q_2^{-1} \right) \sum_{i=1}^n z_i \right], g \left[ X + \left(1-q_1^{-1} \right)\left(1- q_2^{-1} \right) \sum_{i=1}^n z_i^{-1} \right] \right \rangle
$$
To evaluate the sum above, we will use the following formula, which one can prove using the machinery of \cite{Shuf} (similarly with the way one proves \eqref{eqn:pairing shuffle 1}). \\

\begin{proposition}

Consider the following element of $\CA^- \cong \CS$, $\forall d_1,\dots,d_n \in \BZ$:
$$
S_{d_1,\dots,d_n} = \esym \left[ z_1^{d_1} \dots z_n^{d_n} \right]
$$
This element admits the following decomposition in the basis \eqref{eqn:basis of shuffle}: 
$$
S_{d_1,\dots,d_n} = \sum_{v \text{ convex}^*\text{ of size }n} \frac {P_{-v}}{z_v} \int_{|y_1| = \dots = |y_n| = 1}^{|q_1|,|q_2|>1} \frac {\Upsilon^+(P_v)(y_1,\dots,y_n) y_1^{-d_1}\dots y_n^{-d_n}}{\prod_{1\leq i \neq j \leq n} \zeta \left( \frac {y_i}{y_j} \right)} \prod_{a=1}^n  \frac {dy_a}{2\pi i y_a}
$$

\end{proposition} 

\noindent Therefore, if we write: 
$$
S(w_1,\dots,w_n) = \sum_{d_1,\dots,d_n \in \BZ} S_{d_1,\dots,d_n} w_1^{d_1}\dots w_n^{d_n}
$$
as a formal series of elements of $\CA^- \cong \CS$, then we conclude the following. \\

\begin{theorem}
\label{thm:up}

For any $f,g \in \Lambda$, we have the following identity in $\CA^- \cong \CS$:
\begin{equation}
\label{eqn:second 3} 
^{\langle f|}{R'}^{|g\rangle} = \sum_{n=0}^\infty \frac {(-u^{-1}q^{\frac 12})^{n}}{n!} \int^{|q_1|,|q_2| > 1}_{|z_1| = \dots = |z_n| = 1} S(z_1,\dots,z_n) \prod_{a=1}^n  \frac {dz_a}{2\pi i z_a}
\end{equation}
$$
\left \langle f \left[ X + \left(1-q_1^{-1} \right)\left(1- q_2^{-1} \right) \sum_{i=1}^n z_i \right], g \left[ X + \left(1-q_1^{-1} \right)\left(1- q_2^{-1} \right) \sum_{i=1}^n z_i^{-1} \right] \right \rangle
$$
Meanwhile, $^{\langle f|}{R''}^{|g\rangle}$ is given by \eqref{eqn:second}. \\

\end{theorem}

\noindent We emphasize the fact that for any fixed $f$ and $g$, the $n$-th summand of \eqref{eqn:second 3} is a finite linear combination of elements of $\CA^- \cong \CS$. To see this, we note that the second line of \eqref{eqn:second 3} is a Laurent polynomial in $z_1,\dots,z_n$, hence only finitely many terms of the formal series $S(z_1,\dots,z_n)$ survive the contour integral. \\

\begin{remark} Formula \eqref{eqn:second 3} is reminiscent of \cite{FJM, FT2}, where the weighted trace of $\ddR^{F^\rightarrow_u}$ in the first tensor factor was computed (e.g. Proposition 3.3 of \cite{FJM}). \\ 
	
\end{remark} 

\subsection{} In \cite{FJMM}, the authors studied the object $^{\langle 1|}{\ddR}^{|g\rangle}$ for an arbitrary $g \in \Lambda$, but instead of regarding it as an element of $\CA^- \cong \CS$, they regard it as an element of:
$$
\CA = \CA^\uparrow \otimes \BQ(q_1,q_2)\Big \langle P_{n,0}, c_1^{\pm 1}, c_2^{\pm 1} \Big \rangle_{n \neq 0} \otimes \CA^\downarrow
$$
where $\CA^\uparrow$ (resp. $\CA^\downarrow$) is the subalgebra generated by $P_{n,m}$ for all $n \in \BZ$ and $m > 0$ (resp. $m<0$). Because of the automorphisms \eqref{eqn:slz}, these subalgebras are also isomorphic to the shuffle algebra and its opposite, so we have an isomorphism:
\begin{equation}
\label{eqn:embedding}
\CA \cong \CS \otimes \BQ(q_1,q_2)\Big \langle P_{n,0}, c_1^{\pm 1}, c_2^{\pm 1} \Big \rangle_{n \neq 0} \otimes \CS^{\text{op}}
\end{equation}
If we compose the inclusion $\CA^- \subset \CA$ with \eqref{eqn:embedding}, then we obtain:
$$
\CA^- \stackrel{\iota}\hookrightarrow \CS \otimes \BQ(q_1,q_2)\Big \langle P_{n,0}, c_1^{\pm 1}, c_2^{\pm 1} \Big \rangle_{n \neq 0} \otimes \CS^{\text{op}}
$$
As a consequence of \eqref{eqn:r tor}, it is not hard to see that:
\begin{equation}
\label{eqn:bethe}
\iota \left( ^{\langle 1|}{\ddR}^{|g\rangle} \right) \in \CS \otimes \BQ(q_1,q_2)\left[ P_{n,0} \right]_{n < 0}
\end{equation}
The element \eqref{eqn:bethe} is related to off-shell Bethe vectors in \loccit We do not have a closed formula for it, but we will now explain how formula \eqref{eqn:r tor} allows one to obtain an explicit sum over convex paths. Explicitly, for any $f$ and $g$, we have:
$$
^{\langle f |}\ddR^{|g \rangle} =  \sum_{v = \left\{ \frac {m_1}{n_1} \leq \dots \leq \frac {m_t}{n_t} \right\}, (n_i,m_i) \in \BZ_+^2}  \frac {P_{-v}}{z_v} \cdot \langle f | P_{v} | g \rangle 
$$
Let $f = \bar{\lambda} := \Psi ( |\lambda \rangle)$ and $g = \bar{\mu} := \Psi ( |\mu \rangle)$ be the modified Macdonald polynomials associated to partitions $\lambda$ and $\mu$. Then we can invoke \eqref{eqn:iso actions} to obtain the following:
\begin{equation}
\label{eqn:second 2}
^{\langle \bar{\lambda} |}\ddR^{|\bar{\mu} \rangle} = \sum_{v = \left\{ \frac {m_1}{n_1} \leq \dots \leq \frac {m_t}{n_t} \right\}, (n_i,m_i) \in \BZ_+^2}  \frac {P_{-v}}{z_v} \cdot \langle \lambda | P_{m_1,-n_1}\dots P_{m_t,-n_t} |\mu \rangle q^{\sum_i \frac {n_i \e_{m_i,-n_i}}2}
\end{equation}
where now the matrix coefficients are calculated in the representation $F_u^\uparrow$, i.e. according to formulas \eqref{eqn:action right 1}-\eqref{eqn:action right 4}. Note that all but finitely many convex paths have trivial matrix coefficient for any given $\lambda$ and $\mu$, as can be seen from the fact that $P_{m_t,-n_t}|\mu\rangle$ is a linear combination over Young diagrams with $m_t$ boxes fewer than $\mu$. Akin to \eqref{eqn:r tor}, we have the following factorization of the operator \eqref{eqn:second 2}: \\

\begin{proposition}
\label{prop:up}
	
We have the following formula for the matrix coefficients of $\ddR$ in terms of the decomposition \eqref{eqn:embedding}:
$$
^{\langle \bar{\lambda} |}\ddR^{|\bar{\mu} \rangle} = \sum_{\text{Young diagrams }\nu \subset \lambda, \mu} A_{\lambda,\nu} \cdot D_\nu \cdot B_{\nu, \mu}
$$
where:
\begin{align}
&A_{\lambda, \nu} = \sum_{v = \left\{\frac {m_1}{n_1} \leq \dots \leq \frac {m_t}{n_t} \right\}, m_i < 0, n_i > 0}  \frac {P_{-v}}{z_v} \cdot \langle \lambda | P_{m_1,-n_1}\dots P_{m_t,-n_t} |\nu \rangle	\label{eqn:l} \\
&B_{\nu, \mu} = \sum_{v = \left\{\frac {m_1}{n_1} \leq \dots \leq \frac {m_t}{n_t} \right\},m_i > 0, n_i \geq 0}  \frac {P_{-v}}{z_v} \cdot \langle \nu | P_{m_1,-n_1}\dots P_{m_t,-n_t} |\mu \rangle q^{-\sum_i \frac {n_i}2}\label{eqn:u}
\end{align}
and $D_\nu = \sum_{\bar{n} = (n_1 \geq \dots \geq n_t)} \frac {P_{-\bar{n},0}}{z_{\bar{n}}} \langle \nu | P_{0,-\bar{n}} | \nu \rangle$. \\
	
\end{proposition} 

\noindent By analogy with Claim \ref{claim:first}, it is easy to see that:
$$
D_\nu = \exp \left[ \sum_{k=1}^\infty \frac {P_{-k,0}}k \cdot q^{\frac k2} \left( \sum_{\sq \text{ o.c. of }\nu} \chi_\sq^{-k} - \sum_{\sq \text{ i.c. of }\nu} \chi_\sq^{-k} \right) \right] 
$$
which matches formula (4.8) of \cite{FJMM} for the operator $L_{\emptyset, \emptyset}$ in the terminology of \loccit As for the operators \eqref{eqn:l} and \eqref{eqn:u}, formulas \eqref{eqn:action right 3} and \eqref{eqn:action right 4} imply that:
\begin{equation}
\label{eqn:l explicit}
A_{\lambda, \nu} = \sum_{\left\{\frac {m_1}{n_1} \leq \dots \leq \frac {m_t}{n_t} \right\}, m_i < 0, n_i > 0} \frac {P_{-n_1,-m_1}\dots P_{-n_t,-m_t}}{z_v} \cdot
\end{equation}
$$
\Upsilon^-(P_{m_1,-n_1} \dots P_{m_t,-n_t}) (\lambda \backslash \nu) \sigma^{|\lanu|} \prod_{\bsq \in \lanu} \frac {\prod_{\sq \text{ o.c. of }\lambda} \left(1 - \frac {\chi_\sq}{\chi_\bsq} \right)}{\prod_{\sq \text{ i.c. of }\lambda} \left(1 - \frac {\chi_\sq}{\chi_\bsq} \right)} 
$$
\begin{equation}
\label{eqn:u explicit}
B_{\nu, \mu} =  \sum_{\left\{\frac {m_1}{n_1} \leq \dots \leq \frac {m_t}{n_t} \right\},m_i > 0, n_i \geq 0} \frac {P_{-n_1,-m_1}\dots P_{-n_t,-m_t}}{z_v} \cdot
\end{equation}
$$
\Upsilon^+(P_{m_1,-n_1} \dots P_{m_t,-n_t}) (\mu \backslash \nu) \bar{\sigma}^{|\munu|} \prod_{\bsq \in \munu} \frac {\prod_{\sq \text{ i.c. of }\nu} \left(1 - \frac {\chi_\sq}{q\chi_\bsq} \right)}{\prod_{\sq \text{ o.c. of }\nu} \left(1 - \frac {\chi_\sq}{q\chi_\bsq} \right)} q^{-\sum_i \frac {n_i}2}
$$
where $z_v$ denotes \eqref{eqn:zv} for the path $v=\{(n_1,m_1),\dots,(n_t,m_t)\}$. Then one can plug in formula \eqref{eqn:explicit pnm} in order to explicitly compute the second lines of \eqref{eqn:l explicit} and \eqref{eqn:u explicit}. \\

\end{document}